\documentclass{amsart}

\newcommand{\ID}{\mathbb D}
\newcommand{\IR}{\mathbb R}
\newcommand{\IC}{\mathbb C}
\newcommand{\IN}{\mathbb N}
\newcommand{\pr}{\mathrm{pr}}
\newcommand{\Iso}{\mathrm{Iso}}

\newcommand{\non}{\mathrm{non}}
\newcommand{\M}{\mathcal M}
\newtheorem{theorem}{Theorem}
\newtheorem{problem}{Problem}
\newtheorem{proposition}{Proposition}
\newtheorem{lemma}{Lemma}
\theoremstyle{definition}
\newtheorem{remark}{Remark}

\title[On colorings of the Lobachevsky plane]{Symmetric monochromatic subsets in colorings of the Lobachevsky plane}
\author{Taras Banakh}
\address[T.~Banakh]{Department of Mathematics,
 Lviv National University, Lviv, Ukraine,\newline  andInstytut Matematyki, Akademia \'Swi\c etokrzyska, Kielce, Poland}
 \email{tbanakh@yahoo.com}
\author{Artem Dudko}
\address[A.~Dudko]{Kharkiv National University,
 Kharkiv, Ukraine}
\email{artemdudko@rambler.ru}
\author{Du\v san Repov\v s}
\address[D.~Repov\v s]{Institute for Mathematics, Physics and Mechanics, University of Ljubljana, Jadranska 19, Ljubljana, Slovenia}
\email{dusan.repovs@fmf.uni-lj.si}

\begin{document}
\begin{abstract}
We prove that for each partition of the Lobachevsky plane into finitely many Borel pieces one of the cells of the partition contains an unbounded centrally symmetric subset.
\end{abstract}

\maketitle

It follows from \cite{B} (see also \cite[Theorem 1]{BP1}) that for each partition of the $n$-dimensional space $\IR^n$
into $n$ pieces one of the pieces contains an unbounded centrally symmetric subset. On the other hand, $\IR^n$ admits a partition into $(n+1)$ Borel pieces containing no unbounded centrally symmetric subset. For $n=2$ such a partition is drawn at the picture:

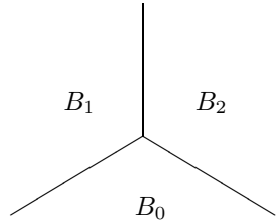
\begin{picture}(300,120)(-180,-50)
\put(0,0){\line(0,1){50}}
\put(0,0){\line(5,-3){50}}
\put(0,0){\line(-5,-3){50}}
\put(-3,-30){$B_0$}
\put(-30,10){$B_1$}
\put(20,10){$B_2$}

\end{picture}

Taking the same partition of the Lobachevsky plane $H^2$, we can see that each cell $B_i$ does contain a unbounded centrally symmetric subset (for such a set just take any hyperbolic line lying in $B_i$). We call a subset $S$ of the hyperbolic plane $H^2$ {\em centrally symmetric or else symmetric with respect to a point} $c\in H^2$ if $S=f_c(S)$ where $f_c:H^2\to H^2$ is the involutive isometry of $H^2$ assigning to each point $x\in H^2$ the unique point $y\in H^2$  such that $c$ is the midpoint of the segment $[x,y]$. The map $f_c$ is called the {\em central symmetry} of $H^2$ with respect to the point $c$.

The following theorem shows that the Lobachevsky plane differs dramatically from the Euclidean plane from the Ramsey point of view.

\begin{theorem}\label{main} For any partition $H^2=B_1\cup\dots\cup B_m$ of the Lobachevsky plane into finitely many Borel pieces one of the pieces contains an unbounded centrally symmetric subset.
\end{theorem}

\begin{proof} We shall prove a bit more: given a partition $H^2=B_1\cup\dots\cup B_m$
of the Lobachevsky plane into $m$ Borel pieces we shall find $i\le m$ and an unbounded subset $S\subset B_i$ symmetric with respect to some point $c$ in an arbitrarily small neighborhood of some finite set $F\subset H^2$ depending only on $m$.

To define this set $F$ it will be convenient to work in the Poincar\'e model of the Lobachevsky plane $H^2$. In this model the hyperbolic plane $H^2$ is identified with the unit disk $\ID=\{z\in\IC:|z|<1\}$ on the complex plane and hyperbolic lines are just segments of circles orthogonal to the boundary of $\ID$.
Let $\overline{\ID}=\{z\in\IC:|z|\le1\}$ be the hyperbolic plane $\ID$ with attached ideal line. For a real number $R>0$ the set $\ID_R=\{z\in\IC:|z|\le 1-1/R\}$ can be thought as a hyperbolic disk of increasing radius as $R$ tends to $\infty$.

On the boundary of the unit disk $\ID$ consider the $(m+1)$-element set
$$A=\{z\in \IC:z^{m+1}=1\}.$$ For any two distinct points $x,y\in A$ by $[x|y]\in\ID$ we denote the ``Euclidean'' midpoint of the arc in $\overline{\ID}$ that connects the points $x,y$ and lies on a hyperbolic line in $H^2=\ID$.
Then $F=\{[x|y]:x,y\in A, \; x\ne y\}$ is a finite subset of cardinality $|F|\le m(m+1)/2$ in the unit disk $\ID$.

We claim that for any open neighborhood $W$ of $F$ in $\IC$ one of the cells of the partition $H^2=B_1\cup\dots\cup B_m$ contains an unbounded subset symmetric with respect to some point $c\in W$. To derive a contradiction we assume the converse and conclude that for every point $c\in W$ and every $i\le m$ the set $B_i\cap f_c(B_i)$ is bounded in $H^2$.

For every $n\in\IN$ consider the set
$$C_n=\{c\in W:\bigcup_{i=1}^m B_i\cap f_c(B_i)\subset \ID_n\}.$$

We claim that $C_n$ is a coanalytic subset of $W$. The latter means that the complement $W\setminus C_n$ is analytic, i.e., is the continuous image of a Polish space. Observe that
$$W\setminus C_n=\{c\in W:\exists i\le m\;\exists x\in \ID\setminus \ID_n,\; x\in B_i \mbox{ and }x\in f_c(B_i)\}=\pr_2(E)$$ where $\pr_2:\ID\times \ID\to \ID$ is the projection on the second factor and
$$E=\bigcup_{i=1}^m\{(x,c)\in \ID\times W: x\in \ID\setminus \ID_n,\; x\in B_i \mbox{ and }f_c(x)\in B_i\}$$is a Borel subset of $\ID\times W$. Being a Borel subset of the Polish space $\ID\times W$, the space $E$ is analytic and so is its continuous image $\pr_2(E)=W\setminus C_n$. Then $C_n$ is coanalytic and hence has the Baire property \cite[21.6]{Ke}, which means that $C_n$ coincides with an open subset $U_n$ of $W$ modulo some meager set. The latter means that the symmetric difference $U_n\triangle C_n$ is meager (that is, of the first Baire category in $W$). Since $C_n\subset C_{n+1}$, we may assume that $U_n\subset U_{n+1}$ for all $n\in\IN$. Let $U=\bigcup_{n=1}^\infty U_n$ and $M=\bigcup_{n=1}^\infty U_n\triangle C_n$.

Taking into account that $W=\bigcup_{n=1}^\infty C_n$, we conclude that
$$W\setminus U=\bigcup_{n=1}^\infty C_n\setminus \bigcup_{n=1}^\infty U_n\subset \bigcup_{n=1}^\infty C_n\setminus U_n\subset \bigcup_{n=1}^\infty C_n\triangle U_n=M$$ which implies that the open set $U$ has meager complement and thus is dense in $W$.

We claim that $F\subset h^{-1}(U)$ for some isometry $h$ of the hyperbolic plane $H^2=\ID$.

For this consider the natural action $$\mu:\Iso(H^2)\times \ID\to \ID,\;\;\mu:(h,x)\mapsto h(x)$$ of the isometry group $\Iso(H^2)$ of the hyperbolic plane $H^2=\ID$. It is easy to see that for every $x\in \ID$ the map $\mu_x:\Iso(H^2)\to \ID$, $\mu_x:h\mapsto h(x)$, is continuous and open (with respect to the compact-open topology on $\Iso(H^2)$). It follows that the set $$\bigcap_{x\in F}\mu_x^{-1}(W)=\{h\in\Iso(H^2):f(F)\subset W\}$$ is an open neighborhood of the neutral element of the group $\Iso(H^2)$.

Taking into account that $U$ is open and dense in $W$ and for every $x\in F$ the map $\mu_x:\Iso(H^2)\to \ID$ is open, we conclude that the preimage
the set $\mu_x^{-1}(U)$ is open and dense in $\mu_x^{-1}(W)\subset\Iso(H^2)$. Then the intersection $\bigcap_{x\in F}\mu_x^{-1}(U)$, being an open dense subset of $\bigcap_{x\in F}\mu^{-1}_x(W)$, is not empty and hence contains some isometry $h$ having the desired property: $F\subset h^{-1}(U)$. Since $F$ is finite, there is $R\in\IN$ with $F\subset h^{-1}(U_R)$. For a complex number $r\in\ID$  consider the set $rA=\{rz:z\in A\}\subset \ID$ and let $$F_r=\{[x|y]:x,y\in rA,\;x\ne y\}\subset\ID,$$ where $[x|y]$ stands for the midpoint of the hyperbolic segment connecting $x$ and $y$ in $H^2$.  It can be shown that for any distinct points
$x,y\in A$ the ``hyperbolic" midpoint $[rx|ry]$ tends to the ``Euclidean" midpoint $[x|y]$ as $r$ tends to 1. Such a continuity yields a neighborhood $O_1$ of 1 such that $F_r\subset h^{-1}(U_R)$ for all $r\in O_1\cap\ID$.

It is clear that for any points $x,y\in A$ the map $$f_{x,y}:\ID\to\ID,\;f_{x,y}:r\mapsto [rx|ry]$$is open and continuous. Consequently, the preimage $f^{-1}_{x,y}(h^{-1}(M))$ is a meager subset of $\ID$ and so is the union $M'=\bigcup_{x,y\in A}f^{-1}_{x,y}(h^{-1}(M))$. So, we can find a non-zero point $r\in O_1\setminus M'$ so close to 1 that the set $rA$ is disjoint with the hyperbolic disk $h^{-1}(\ID_R)$. For this point $r$ we shall get $F_r\cap h^{-1}(M)=\emptyset$.

The set $rA$ consists of $m+1$ points. Consequently, some cell $h^{-1}(B_i)$ of the partition $\ID=h^{-1}(B_1)\cup\dots\cup h^{-1}(B_m)$ contains two distinct points $rx,ry$ of $rA$. Those points are symmetric with respect to the point $$[rx|ry]\in F_r\subset h^{-1}(U_R)\setminus h^{-1}(M).$$ Then the images $a=h(rx)$ and $b=h(ry)$ belong to $B_i$ and are symmetric with respect to the point $c=h([rx|ry])\in U_R\setminus M\subset C_R$. It follows from the definition of $C_R$ that $\{a,b\}\subset B_i\cap f_c(B_i)\subset \ID_R$, which is not the case because $rx,ry\notin h^{-1}(\ID_R)$.
\end{proof}

We do not know if Theorem~\ref{main} is true for any finite (not necessarily Borel) partition of the Lobachevsky plane $H^2$. For partitions of $H^2$ into two pieces the Borel assumption is supefluous.

\begin{theorem}\label{main2} There is a subset $T\subset H^2$ of cardinality $|T|=3$ such that for any partition $H^2=A_1\cup A_2$ of $H^2$ into two pieces either $A_1$ or $A_2$ contains an unbounded subset, symmetric with respect to some point $c\in T$.
\end{theorem}

\begin{proof} Lemma~\ref{lem} below allows us to find an equilateral triangle $\triangle c_0c_1c_2$ on the Lobachevsky plane $H^2$ such that the composition $f_{c_2}\circ f_{c_1}\circ f_{c_0}$ of the symmetries with respect to the points $c_0,c_1,c_2$ coincides with the rotation on the angle $2\pi/3$ around some point $o\in H^2$.
Consequently $(f_{c_2}\circ f_{c_1}\circ f_{c_0})^3$ is the identity isometry of $H^2$.

We claim that for any partition $H^2=A_1\sqcup A_2$ of the Lobachevsky plane into two pieces one of the pieces contains an unbounded subset symmetric with respect to some point in the triangle $T=\{c_0,c_1,c_2\}$. Assuming the converse, we conclude that the set
$$B=\bigcup_{c\in T}\bigcup_{i=1}^2A_i\cap f_c(A_i)$$
is bounded. It follows that two points $x,y\in H^2\setminus B$, symmetric with respect to a center $c\in T$ cannot belong to the same cell $A_i$ of the partition.

Given a point $x_0\in H^2$ consider the sequence of points $x_1,\dots x_9$ defined by the recursive formula:
$x_{i+1}=f_{c_{i\,\mathrm{mod}\,3}}(x_i)$.
It follows that $x_9=(f_{c_2}\circ f_{c_1}\circ f_{c_0})^3(x_0)=x_0$. Taking $x_0$ sufficiently far from the center $o$ of rotation we can guarantee that none of the points $x_0,\dots,x_9$ belongs to $B$.

The point $x_0$ belongs either to $A_1$ or to $A_2$. We lose no generality assuming that $x_0\in A_2$. Since the points $x_0,x_1\notin B$ are symmetric with respect to $c_0$ and $x_0\in A_2$, we get that $x_1\in H^2\setminus A_2=A_1$. By the same reason $x_1,x_2$ cannot simultaneously belong to $A_1$ and hence $x_2\in A_2$. Continuing in this fashion we conclude that $x_i$ belongs to $A_1$ for odd $i$ and to $A_2$ for even $i$. In particular, $x_9\in A_1$, which is not possible because $x_9=x_0\in A_2$.
\end{proof}

\begin{lemma}\label{lem} There is an equilateral triangle $\triangle ABC$ on the Lobachevsky plane such that  the composition $f_{C}\circ f_{B}\circ f_{A}$ of the symmetries with respect to the points $A,B,C$ coincides with the rotation on the angle $2\pi/3$ around some point $O$.
\end{lemma}

\begin{proof} For a positive real number $t$ consider an equilateral triangle $\triangle ABC$ with side $t$ the on the Lobachevsky plane. Let $M$ be the midpoint of the side $AB$ and $l$ be the line through $C$ that is orthogonal to the line $CM$.
Consider also the line $p$ that is orthogonal to the line $AB$ and passes through the point $P$  such that $A$ is the midpoint between $P$ and $M$. Observe that $|PM|=|AB|=t$ and for sufficiently small $t$ the lines $p$ and $l$ intersect at some point $O$.

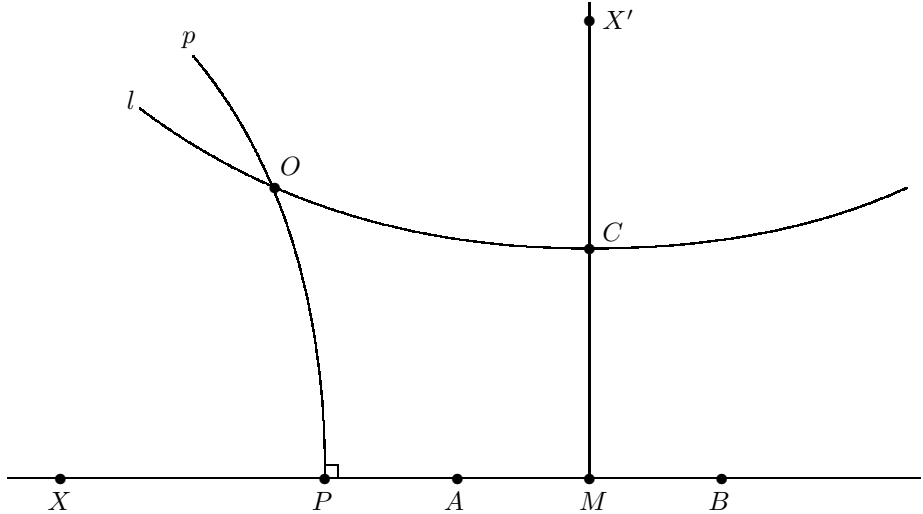
\begin{picture}(300,220)(-220,-20)
\put(0,87){\circle*{4}}
\put(5,90){$C$}
\qbezier(0,87)(70,87)(120,110)
\qbezier(0,87)(-100,87)(-170,140)
\put(-175,140){$l$}
\put(-50,0){\circle*{4}}
\put(-55,-12){$A$}
\put(50,0){\circle*{4}}
\put(45,-12){$B$}
\put(-220,0){\line(1,0){345}}

\put(0,0){\circle*{4}}
\put(-4,-12){$M$}
\put(0,173){\circle*{4}}
\put(5,170){$X'$}
\put(-100,0){\circle*{4}}
\put(-105,-12){$P$}
\qbezier(-100,0)(-100,100)(-150,160)
\put(-95,0){\line(0,1){5}}
\put(-95,5){\line(-1,0){5}}

\put(-154,165){$p$}
\put(-200,0){\circle*{4}}
\put(-205,-12){$X$}
\put(0,0){\line(0,1){180}}
\put(-119,110){\circle*{4}}
\put(-117,115){$O$}
\end{picture}

It is easy to see that the composition $f_B\circ f_A$ is the shift along the line $AB$ on the distance $2t$ and hence the image $f_B\circ f_A(O)$ of the point $O$ is the point symmetric to $O$ with respect to the point $C$. Consequently, $f_C\circ f_B\circ f_A(O)=O$, which means that the isometry $f_C\circ f_B\circ f_A$ is a rotation of the Lobachevsky plane around the point $O$ on some angle $\varphi_t$.

 To estimate this angle, consider the point $X$ such that $P$ is the midpoint between $X$ and $M$. Then $|XM|=2t$ and consequently, $f_B\circ f_A(X)=M$ while $X'=f_C\circ f_B\circ f_A=f_C(M)$ is the point on the line $CM$ such that $C$ is the midpoint between $X'$ and $M$. It follows that $|X'X|\le |XM|+|MX'|<2t+2t=4t$.

Observe that for small $t$ the point $X'$ is near to the point, symmetric to $X$ with respect to $O$, which means that the angle $\varphi_t=\angle XOX'$ is close to $\pi$ for $t$ close to zero.
On the other hand, for very large $t$ the lines $p$ and $l$ on the Lobachevsky plane do not intersect. So we can consider the smallest upper bound $t_0$ of numbers $t$ for which the lines $l$ and $p$ meet. For values $t<t_0$ near to $t_0$  the point $O$  tends to infinity as $t$ tends to $t_0$. Since the length of the side $XX'$ of the triangle $\triangle XOX'$ is bounded by $4t_0$ the angle $\varphi_t=\angle XOX'$ tends to zero as $O$ tends to infinity. Since the angle $\varphi_t$ depends continuous on $t$ and decreases from $\pi$ to zero as $t$ increases from zero to $t_0$, there is a value $t$ such that $\varphi_t=2\pi/3$. For such $t$ the composition $f_C\circ f_B\circ f_A$ is the rotation around $O$ on the angle $2\pi/3$.
\end{proof}

\section*{Some comments and Open Problems}

In contrast to Theorem~\ref{main}, Theorem~\ref{main2} is true for the Euclidean plane $E^2$ even in a stronger form: for any subset $C\subset E^2$ not lying on a line and any partition $E^2=A_1\cup A_2$ one of the cells of the partition contains an unbounded subset symmetric with respect to some center $c\in C$, see \cite{B2}.

Having in mind this result let us call a subset $C$ of a Lobachevsky or Euclidean space $X$ {\em central for (Borel) $k$-partitions} if for any partition $X=A_1\cup\dots\cup A_k$ of $X$ into $k$ (Borel) pieces one of the pieces contains an unbounded monochromatic subset $S\subset X$, symmetric with respect to some point $c\in C$. By $c_k(X)$ (resp. $c^B_k(X)$) we shall denote the smallest size of a  subset $C\subset X$, central for (Borel) $k$-partitions of $X$. If no such a set $C$ exists, then we put $c_k(X)=\infty$ (resp. $c_k^B(X)=\infty$) where $\infty$ is assumed to be greater than any cardinal number. It follows from the definition that  $c_k^B(X)\le c_k(X)$.

We have a lot of information on the numbers $c_k^B(E^n)$ and $c_k(E^n)$ for Euclidean spaces $E^n$, see \cite{B2}. In particular, we known that
\begin{enumerate}
\item  $c_2(E^n)=c_2^B(E^n)=3$ for all $n\ge 2$;
\item $c_3(E^3)=c_3^B(E^3)=6$;
\item $12\le c_4^B(E^4)\le c_4(E^4)\le 14$;
\item $n(n+1)/2\le c_n^B(E^n)\le c_n(E^n)\le 2^n-2$ for every $n\ge 3$.
\end{enumerate}

Much less is known on the numbers $c_k^B(H^n)$ and $c_k^B(H^n)$ in the hyperbolic case.
Theorem~\ref{main2} yields the upper bound $c_2(H^2)\le 3$. In fact, 3 is the exact value of $c_2(H^n)$ for all $n\ge 2$.

\begin{proposition} $c_2^B(H^n)=c_2(H^n)=3$ for all $n\ge 2$.
\end{proposition}

\begin{proof} The upper bound $c_2(H^n)\le c_2(H^2)\le 3$ follows from Theorem~\ref{main2}. The lower bound $3\le c_2^B(H^n)$ will follow as soon as for any two points $c_1,c_2\in H^n$ we construct a partition $H^n=A_1\cup A_2$ in two Borel  pieces containing no unbounded set, symmetric with respect to a point $c_i$. To construct such a partition, consider the line $l$ containing the points $c_1,c_2$ and decompose $l$ into two half-lines $l=l_1\sqcup l_2$. Next, let $H$ be an $(n-1)$-hyperplane in $H^n$, orthogonal to the line $l$. Let $S$ be the unit sphere in $H$ centered at the intersection point of $l$ and $H$. Let $S=B_1\cup B_2$ be a partition of $S$ into two Borel pieces such that no antipodal points of $S$ lie in the same cell of the partition.
For each point $x\in H^n\setminus l$ consider the hyperbolic plane
$P_x$ containing the points $x,c_1,c_2$. The complement
$P_x\setminus l$ decomposes into two half-planes $P^+_x\cup P^-_x$
where $P_x^+$ is the half-plane containing the point $x$. The plane
$P_x$ intersects the hyperplane $H$ by a hyperbolic line containing
two points of the sphere $S$. Finally put
$$A_i=l_i\cup\{x\in H^2\setminus l: P^+_x\cap B_i\ne\emptyset\}$$for $i\in\{1,2\}$.
It is easy to check that $A_1\sqcup A_2=H^n$ is the desired partition of the hyperbolic space into two Borel pieces none of which contains an unbounded subset symmetric with respect to one of the points $c_1$, $c_2$.
\end{proof}

The preceding proposition implies that the cardinal numbers $c_2(H^n)$ are finite.

\begin{problem} For which numbers $k,n$ are the cardinal numbers $c_k(H^n)$ and $c_k^B(H^n)$ finite? Is it true for all $k\le n$?
\end{problem}

Except for the equality $c_2(E^n)=3$, we have no information on the numbers $c_k(E^n)$ with $k<n$.

\begin{problem} Calculate (or at least evaluate) the numbers $c_k(E^n)$ and $c_k(H^n)$ for $2<k<n$.
\end{problem}

In all the cases where we know the exact values of the numbers $c_k(E^n)$ and $c_k^B(E^n)$ we see that those numbers are equal.

\begin{problem} Are the numbers $c_k(E^n)$ and $c_k^B(E^n)$ (resp. $c_k(H^n)$ and $c_k^B(H^n)$) equal for all $k,n$?
\end{problem}

Having in mind that each subset not lying on a line is central for 2-partitions of the Euclidean plane, we may ask about the same property of the Lobachevsky plane.

\begin{problem} Is any subset $C\subset H^2$ not lying on a line central for (Borel) 2-partitions of the Lobachevsky plane $H^2$?
\end{problem}

Finally, let us ask about the numbers $c_k^B(H^2)$ and $c_k(H^2)$.
Observe that Theorem~\ref{main} guarantees that $c^B_k(H^2)\le \mathfrak c$ for all $k\in\IN$. Inspecting the proof we can see that this upper bound can be improved to $c^B_k(H^2)\le\non(\mathcal M)$ where $\non(\M)$ is the smallest cardinality of a non-meager subset of the real line. It is clear that $\aleph_1\le\non(\mathcal M)\le \mathfrak c$. The exact location of the cardinal $\non(M)$ on the interval $[\aleph_1,\mathfrak c]$ depends on axioms of Set Theory, see \cite{Bl}. In particular, the inequality $\aleph_1=\non(\M)<\mathfrak c$ is consistent with ZFC.

\begin{problem} Is the inequality $c_k^B(H^2)\le \aleph_1$ provable in ZFC? Are the cardinals $c_k^B(H^2)$ countable? finite?
\end{problem}

The latter problem asks if $H^2$ contains a countable (or finite) central set for Borel $k$-partitions of the Lobachevsky plane. Inspecting the proof of Theorem~\ref{main} we can see that it gives an ``approximate'' answer to this problem:

\begin{proposition}\label{main3} For any $k\in\IN$ there is a finite subset $C\subset H^2$of cardinality $|C|\le k(k+1)/2$ such that for any partition $H^2=B_1\cup\dots\cup B_k$ of $H^2$ into $k$ Borel pieces and for any open neighborhood $O(C)\subset H^2$ of $C$ one of the pieces $B_i$ contains an unbounded subset $S\subset B_i$ symmetric with respect to some point $c\in O(C)$.
\end{proposition}




\begin{remark} For further results and open problems related to symmetry and colorings see the surveys \cite{BP2}, \cite{BVV} and the list of problems \cite[\S4]{BBGRZ}.
\end{remark}

\end{document}